\documentclass[11pt,twoside,leqno]{amsart}

\usepackage{amscd, amsfonts, amsmath, amssymb, amsthm, bm, epsf, epsfig, flafter, float, subfigure}
\usepackage[all]{xy}
\usepackage[margin=1in]{geometry}

\theoremstyle{plain}
\newtheorem{theorem}{Theorem}[section]

\newtheorem*{claim*}{Claim}

\newtheorem{question}[theorem]{Question}

\theoremstyle{definition}
\newtheorem{definition}[theorem]{Definition}

\theoremstyle{remark}
\newtheorem{remark}[theorem]{Remark}

\newcommand{\rmEr} {{I\kern-1pt I }}
\newcommand{\rmSan}{{I\kern-1pt I\kern-1pt I }}
\newcommand{\rmSi} {{I\kern-2pt V }}

\begin{document}

\title[The support genus of certain Legendrian knots]{The support genus of certain Legendrian knots}

\author[Youlin Li]{Youlin Li}
\address{Department of Mathematics \\Shanghai Jiaotong University\\
Shanghai, 200240, P.R.China}
\email{liyoulin@sjtu.edu.cn}

\author[Jiajun Wang]{Jiajun Wang}
\address{Beijing International Center for Mathematical Research \& School of Mathematical Sciences, Peking University, Beijing, 100871, P.R.China}
\email{wjiajun@math.pku.edu.cn}


\begin{abstract} In this paper, the support genus of all Legendrian right handed trefoil knots and some other Legendrian knots is computed. We give examples of Legendrian knots in the three-sphere with the standard contact structure which have positive support genus with arbitrarily negative Thurston-Benniquin invariant. This answers a question in \cite{o09}.
\end{abstract}

\maketitle

\section{\textbf{Introduction}}

In his seminal paper \cite{g02}, Giroux established the surprising
one-to-one correspondence between the contact structures up to
isotopy and the open book decompositions up to positive
stabilizations on a given closed oriented three-manifold. See
\cite{E06} for details. It becomes natural and convenient for
topologists to study contact structures in the viewpoint of  open
book decompositions.


For Legendrian knots, Akbulut and Ozbagci \cite{AO01} showed that for any Legendrian link $L$ in $S^3$ with the standard contact structure $\xi_{\rm std}$, there exists a compatible open book such that $L$ sits in a page of the open book, furthermore, the framing of L given by the page of the open book agrees with the contact framing. See also \cite{p04}.

In \cite{eo08}, Etnyre and Ozbagci introduced the
definition of support genus for a contact three-manifold, which is
the minimal genus of a page among all open books supporting the
given the contact three-manifold. 
 In a similar fashion, Onaran
defined in \cite{o09} an invariant for Legendrian knots in a contact
three-manifold as follows.

\begin{definition}[\cite{o09}] Let $L$ be a Legendrian knot in a  contact three-manifold $(M,\xi)$, the \emph{support genus} of $L$, denoted by $sg(L)$, is the minimal genus of a page among all open book decompositions of $M$ supporting $\xi$ such that $L$ sits on a page of the open book and the framings given by $\xi$ and given by the page coincide.
\end{definition}

Ding and Geiges \cite{dg04} introduced the definition of contact surgeries along a Legendrian link. If a contact three-manifold $(M, \xi)$ is
obtained by a contact $r$ surgery along a Legendrian knot $L$ in
$(S^{3}, \xi_{\rm std})$, Onaran showed that $sg(L)$ is greater than
or equal to the support genus of $(M, \xi)$ (\cite[Remark
5.11]{o09}). In the same paper, the following question is asked:

\begin{question}
Does every Legendrian knot in $(S^{3}, \xi_{\rm std})$ with negative Thurstion-Bennequin invariant have support genus zero?
\end{question}

In the present paper, we study the support genus of certain Legendrian knots in $(S^{3}, \xi_{\rm std})$. First, for Legendrian torus knots, we have the following

\begin{theorem}
Suppose $k\geq1$. Let $L$ be a Legendrian torus knot $T(2,2k+1)$ in
$(S^{3}, \xi_{\rm std})$ with nonnegative Thurston-Bennequin
invariant, then $sg(L)=1$.
\end{theorem}

Next, we study the support genus of Legendrian twist knots
$K_{-2m}$, where $m\geq1$. (See \cite{env10} for the meaning of
$K_{-2m}$.) In particular, $K_{-2}$ is the right handed trefoil. All
Legendrian twist knots are classified in \cite{env10}. In fact,
$K_{-2m}$ has $\lceil m^{2}/2\rceil$ Legendrian representatives with
Thurston-Bennequin invariant one and rotation number zero, and has a
unique Legendrian representative with Thurston-Bennequin invariant
minus one and rotation number zero (\cite[Theorem 1.1 (4)]{env10}).
We have the following

\begin{theorem}
Suppose $m, n_1$ and $n_2$ are natural numbers. Let $L$ be a Legendrian representative of the twist knot $K_{-2m}$ in $(S^{3}, \xi_{\rm std})$ with Thurston-Bennequin invariant $1$. Then $sg(S_{+}^{n_{1}}S_{-}^{n_{2}}(L))=0$, where $S_+$ and $S_-$ denote the positive and negative stabilizations respectively. 
\end{theorem}

Theorems 1.3 and 1.4 both computed the support genus of some Legendrian right-handed trefoil knots. By the classification results in \cite{eh01}, the remaining Legendrian right-handed trefoil knots are $S_{+}^{n}(L)$ and $S_{-}^{n}(L)$ for $n\geq2$. The support genus of these Legendrian knots is computed in the following theorem, and thus we computed the support genus of all Legendrian right-handed trefoil knots.

\begin{theorem}
Let $L$ be a Legendrian right handed trefoil knot in $(S^{3},\xi_{\rm std})$ with Thurston-Bennequin invariant $1$. Then for any integer $n\geq2$, both $S_{+}^{n}(L)$ and $S_{-}^{n}(L)$ have support genus $1$.
\end{theorem}

This gives a negative answer to Question 1.2. In fact, Legendrian
knots with positive support genus can have arbitrarily negative
Thurston-Bennequin invariant.

\vspace{10pt}

Our strategy is the following. We will construct a fibered link with
a Thurston norm minimizing Seifert surface which contains the
interested knot. By computing the Thurston-Bennequin invariant and
rotation number of the Legendrian realization of this knot and known
classification results, we will be able to determine the Legendrian
knot and get an upper bound for the support genus for it. In the
other direction, we use the Heegaard Floer contact invariant to get
lower bound for the support genus. Combining these results, we will
compute the support genus of our interested Legendrian knots.

\subsection*{Acknowledgement} Part of the work was  done when the first named author was visiting Peking University. He would like to thank the School of Mathematical Sciences in Peking University for their hospitality. He is also partially supported by NSFC grant 11001171.

\section{\textbf{Proof of Theorem 1.3}}

\begin{figure}[htbp]
\center{\includegraphics[width=140pt]{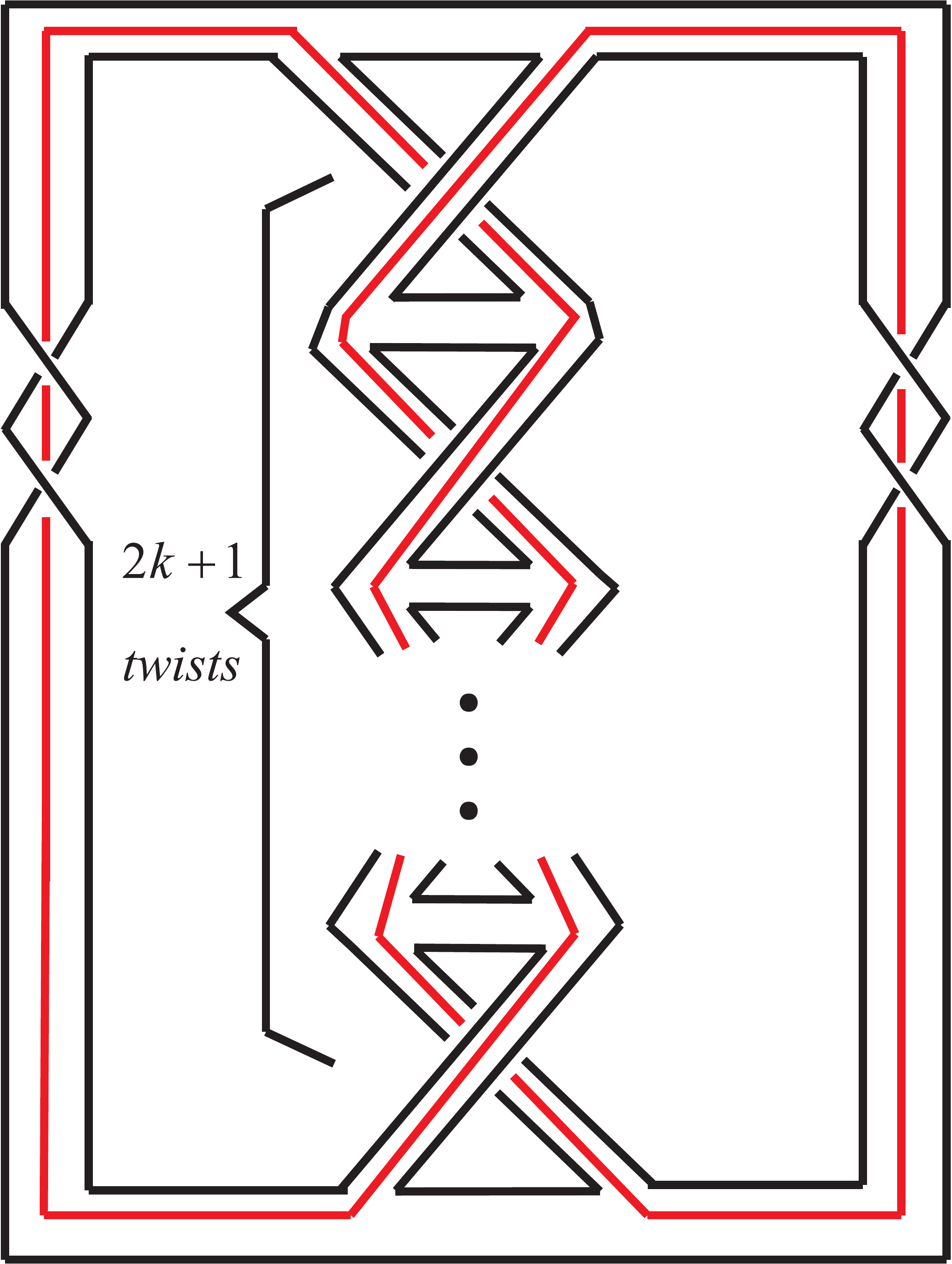}}
\begin{minipage}{6in}
\caption{\label{fig:open_book_2_2kplus1} {\bf A punctured torus
containing the torus knot $T(2,2k+1)$.} }
\end{minipage}
\end{figure}

The surface $F$ illustrated in Figure \ref{fig:open_book_2_2kplus1}
is a punctured torus and contains an embedded torus knot $T(2,2k+1)$
in its interior.

\begin{figure}[htbp]
\center{\includegraphics[width=350pt]{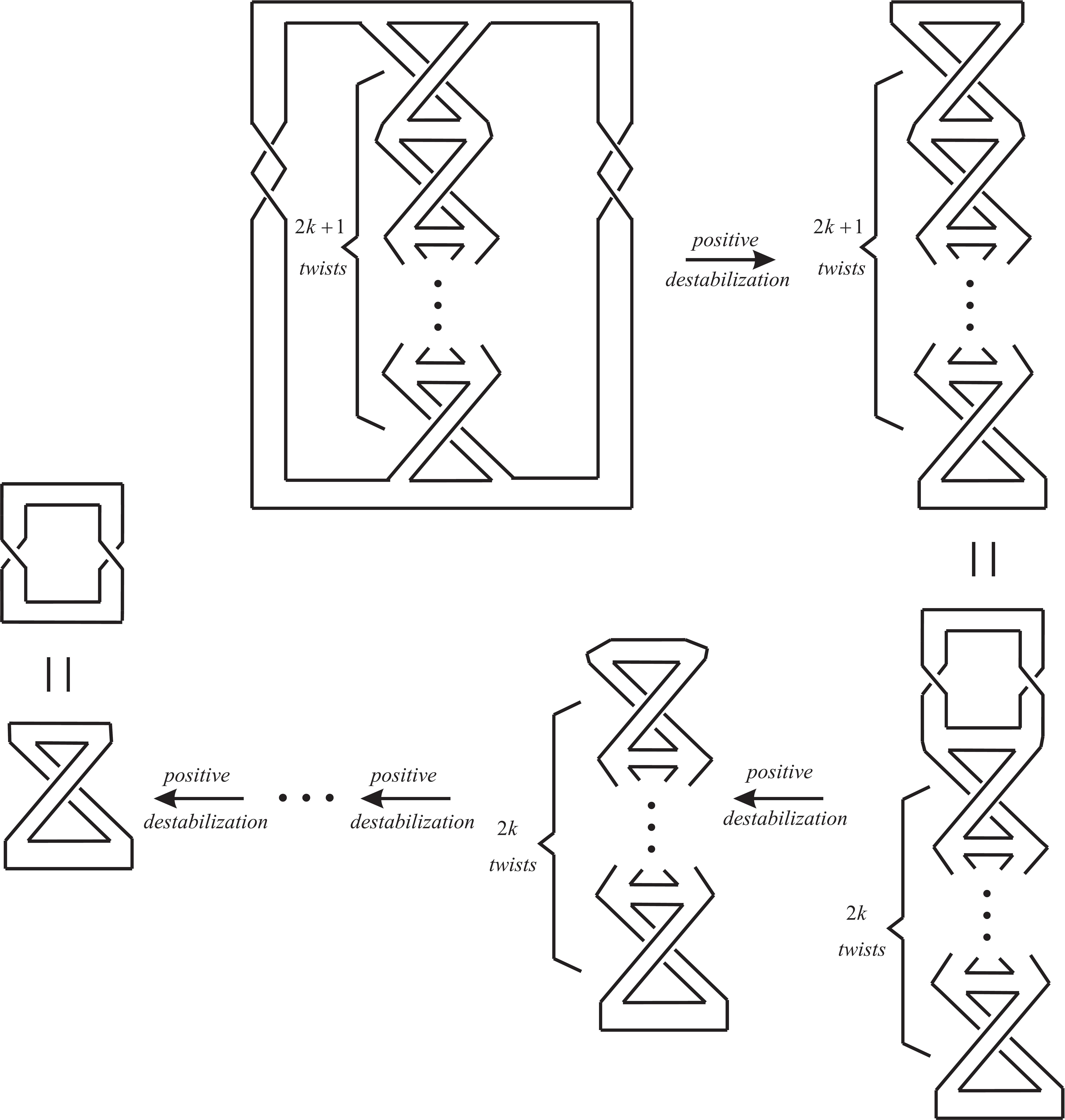}}
\begin{minipage}{6in}
\caption{\label{fig:stabilization_2_2kplus1} {\bf Positive
destabilizations of a page containing the torus knot  $T(2,2k+1)$.}
}
\end{minipage}
\end{figure}

As illustrated in Figure \ref{fig:stabilization_2_2kplus1}, $F$  can
be obtained by a sequence of positive stabilizations from the disk.
By \cite{g02}, $F$ can be viewed as a page of an open book
decomposition of $S^{3}$ which supports the standard tight contact
structure $\xi_{\rm std}$. $K$ can be realized as a Legendrian knot
since $K$ is homologically nontrivial (in fact, nonseparating) in
$F$, which we denote by $T_m(2,2k+1)$. Moreover, the contact framing
of $T_m(2,2k+1)$ coincides with the page framing induced by $F$.
This is because the contact planes can be arranged to be arbitrarily
close to the tangent planes of the pages.

Let $J$ be a push-off of $T_m(2,2k+1)$ along the surface $F$, then
$J$ represents the contact framing of $T_m(2,2k+1)$. The
Thurston-Bennequin invariant of $T_m(2,2k+1)$ equals to the linking
number of $T_m(2,2k+1)$ and $J$, which is easily seen to be $2k-1$.
By the classification of Legendrian representatives of torus knots
in \cite{eh01}, $T_m(2,2k+1)$ has the maximum Thurston-Bennequin
invariant over all Legendrian representatives of $T(2,2k+1)$ and it
is the unique one.

So the support genus of the Legendrian knot $T_m(2,2k+1)$ is at most
$1$. According to \cite{eh01}, all other Legendrian representatives
of the torus knot $T(2,2k+1)$ can be obtained by stabilizing
$T_m(2,2k+1)$. Thus, by \cite[Theorem 5.9]{o09}, the support genus
of any Legendrian representative of $T(2,2k+1)$ is at most $1$.

On the other hand, since a Legendrian knot in a weakly fillable
contact structure with positive Thurston-Bennequin invariant has
positive support genus (\cite[Lemma 5.4]{o09}), the support genus of
$T_m(2,2k+1)$ is at least $1$. Therefore the support genus of
$T_m(2,2k+1)$ is $1$.

This ends the proof of Theorem 1.3.

\section{\textbf{Proof of Theorem 1.4}}

By \cite[Theorem 5.9]{o09}, it suffices to show that the Legendrian
representative of the twist knot $K_{-2m}$ with Thurston-Bennequin
invariant minus one and rotation number zero has support genus zero.

\begin{figure}[htbp]
\center{\includegraphics[width=450pt]{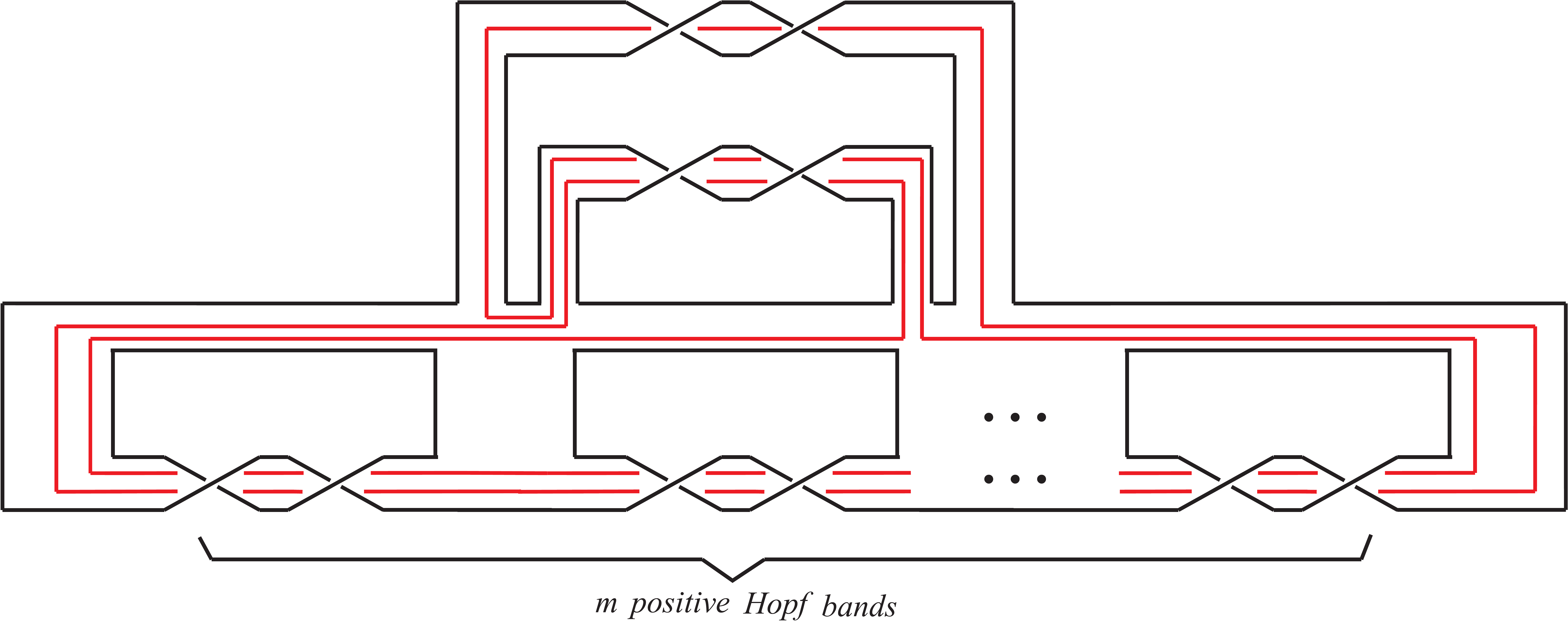}}
\begin{minipage}{6in}
\caption{\label{fig:open_book_twist}
{\bf A punctured sphere containing the twist knot $K_{-2m}$.} }
\end{minipage}
\end{figure}

The punctured sphere illustratd in Figure \ref{fig:open_book_twist}
contains the twist knot $K_{-2m}$ in its interior.
It is easy to see that $F$ is a page of an open book  decomposition
which corresponds to $(S^{3},\xi_{std})$. Since $K$ is not null
homologous in $F$, $K$ can be made Legendrian and the contact
framing of $K$ coincides with the page framing induced by $F$.

The linking number of $K$ and its push-off along $F$ is $-1$,  by
the same argument as in the previous section, the Thurston-Bennequin
invariant of $K$ is $-1$. Below we shall compute the rotation number
of $K$.

\begin{figure}[htbp]
\center{\includegraphics[width=350pt]{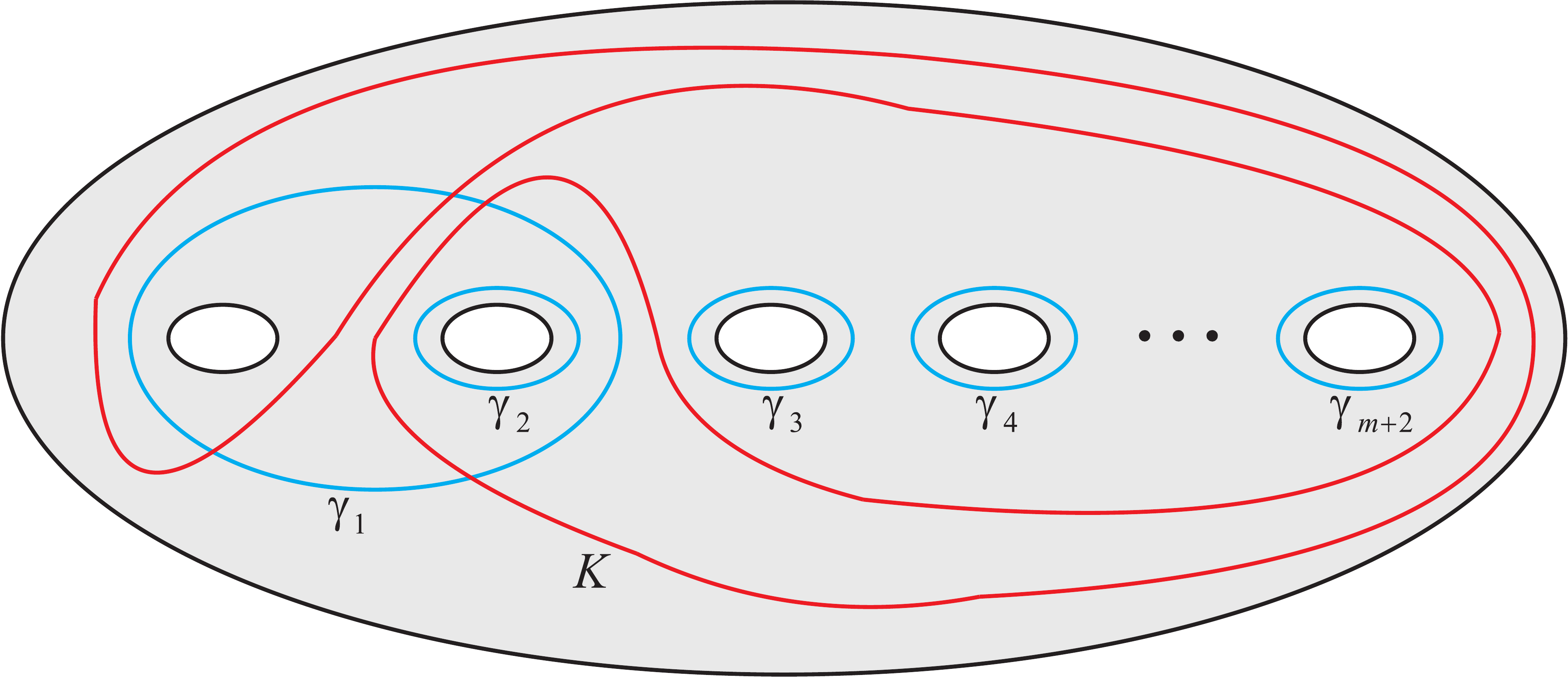}}
\begin{minipage}{6in}
\caption{\label{fig:abstract1_open_book_twist}
{\bf An abstract open book.} }
\end{minipage}
\end{figure}

We turn the open book decomposition shown in Figure
\ref{fig:open_book_twist}  into an abstract open book decomposition
$(F',t_{\gamma_{1}}t_{\gamma_{2}}\ldots t_{\gamma_{m+2}})$ shown in
Figure \ref{fig:abstract1_open_book_twist}, where $F'$ is a
punctured sphere which is homeomorphic to $F$, and $t_{\gamma_{i}}$
denotes the right handed Dehn twist along $\gamma_{i}$, $i=1, 2,
\ldots, m+2$.

If we perform a Legendrian surgery along $K$, then we obtain a Stein
fillable tight contact manifold which corresponds to the planar open
book decomposition $(F', t_{\gamma_{1}}t_{\gamma_{2}}\ldots
t_{\gamma_{m+2}}t_{K})$. Let $(W, J)$ be the Stein surface obtained
from $B^{4}$ by attaching a two-handle corresponding to the Legendrian
surgery along $K$. Let $c_{1}(J)$ be the first Chern class of this
Stein surface, and $h$ be the generator of
$H_{2}(W;\mathbb{Z})\cong\mathbb{Z}$ supported on the attached 2-handle,
then $\left<c_{1}(J), h\right>=rot(K)$.

\begin{figure}[htbp]
\center{\includegraphics[width=350pt]{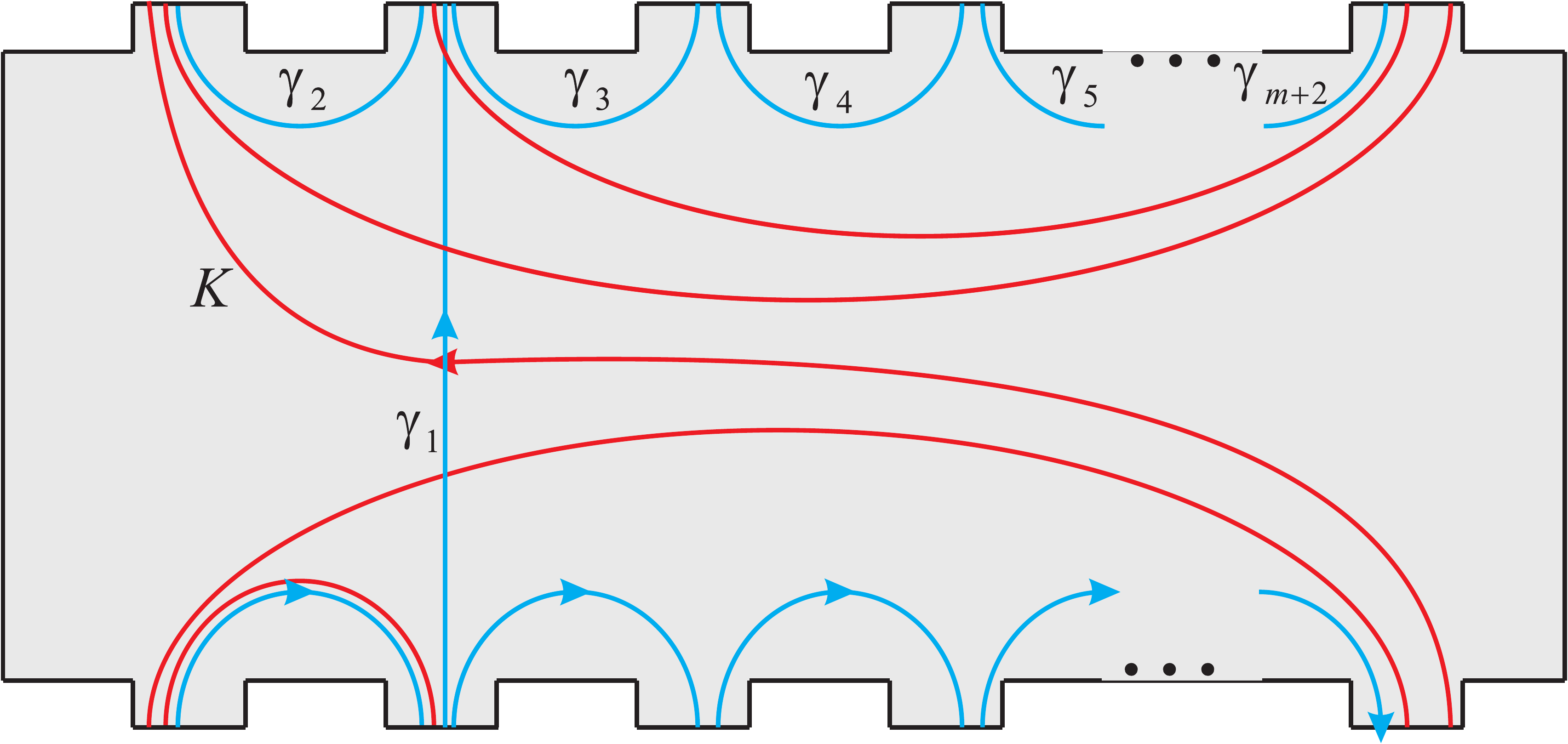}}
\begin{minipage}{6in}
\caption{\label{fig:abstract2_open_book_twist}
{\bf Another diagram for the abstract open book.} }
\end{minipage}
\end{figure}

\begin{figure}[htbp]
\center{\includegraphics[width=320pt]{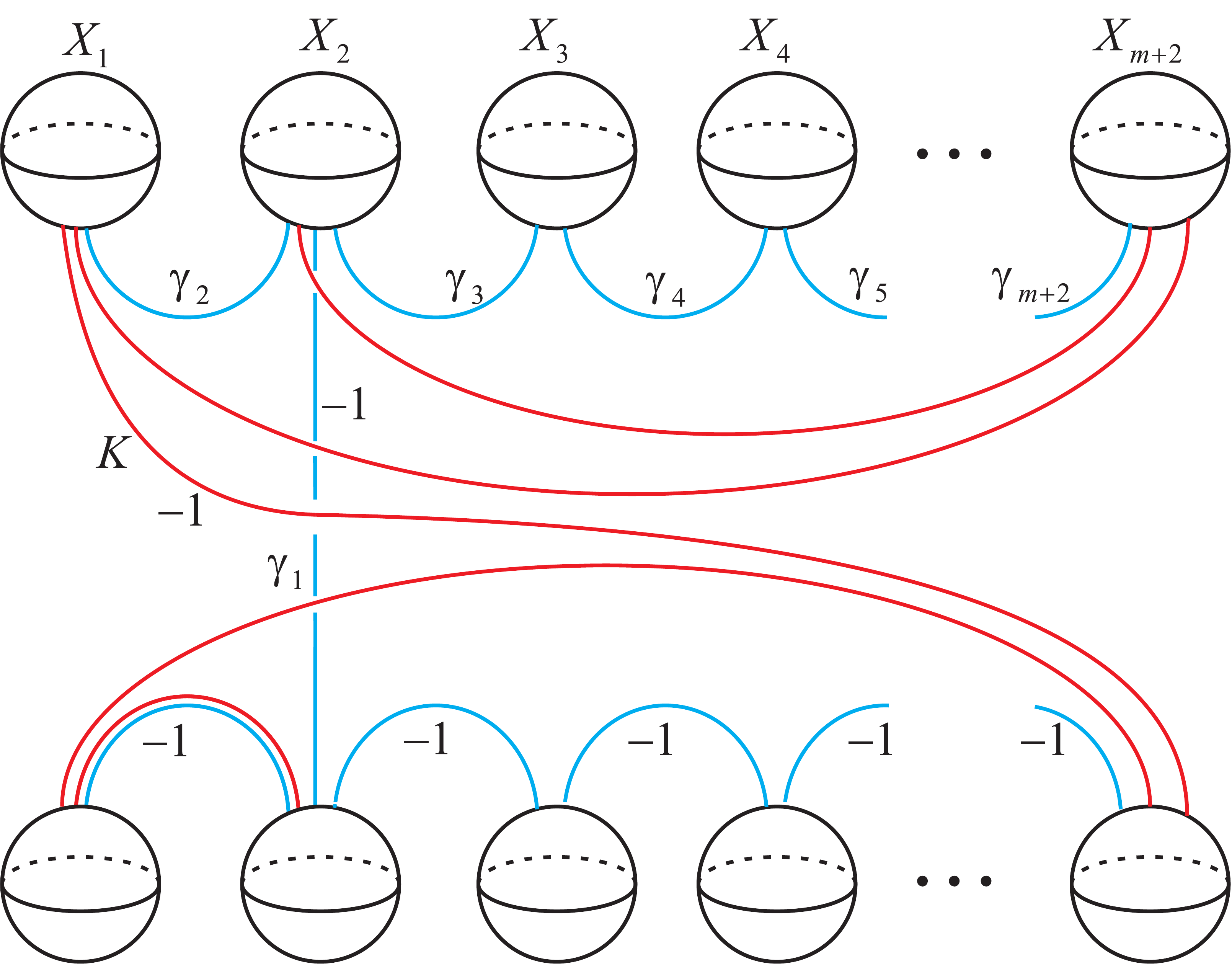}}
\begin{minipage}{6in}
\caption{\label{fig:kirby_diagram_stein}
{\bf A Kirby diagram for the Stein surface.} }
\end{minipage}
\end{figure}

The open book decomposition in Figure
\ref{fig:abstract1_open_book_twist}  is the same as the one in
Figure \ref{fig:abstract2_open_book_twist},
where each upper horizontal segment is identified with the
corresponding  lower horizontal segment to form a one-handle. We
consider the Kirby diagram in Figure \ref{fig:kirby_diagram_stein},
where the framings are labelled with respect to the blackboard
framings.  According to \cite{e90}, this Kirby diagram presents a
Stein surface, denoted by $(W^\prime, J^\prime)$. By handle
cancellation, it is not hard to see that $W$ and $W^\prime$ are
diffeomorphic. Since the induced contact structures on $\partial W$
by $(W, J)$ and $(W^\prime, J')$ are isotopic, we have
$c_{1}(J)=c_{1}(J')$ by \cite[Theorem 1.2]{lm97}.

Let $X_{1}$, $X_{2}$, $\ldots$,  $X_{m+2}$ denote the 1-handles
attached to $D^{4}$ to form $\#^{m+2}(S^{1}\times D^{3})$. Let $S
_{\gamma_{1}}$, $S _{\gamma_{2}}$, $\ldots$, $S _{\gamma_{m+2}}$ and
$S _{K}$ be the cores of the 2-handles attached to the curves
$\gamma_{1}$, $\gamma_{2}$, $\ldots$,  $\gamma_{m+2}$ and $K$,
respectively, and let $C _{\gamma_{1}}$, $C _{\gamma_{2}}$,
$\ldots$, $C_{\gamma_{m+2}}$ and $C_{K}$ be the cocores of the
2-handles.

In the tight contact manifold $\#^{m+2}(S^{1}\times S^{2})$
supported  by the open book decomposition $(F', id)$, we can
Legendrian realize $\gamma_{1}$, $\gamma_{2}$, $\ldots$,
$\gamma_{m+2}$ and $K$. Moreover, by the argument in Section 3.1 in
\cite{eo08} and Figure 5, we can easily compute their rotation
numbers as
$$r(\gamma_{1})=0, r(\gamma_{2})=r(\gamma_{3})=\ldots=r(\gamma_{m+2})=-1$$
and
$r(K)=0$. Also according to Section 3.1 in \cite{eo08}, the class
$c_{1}(J')$ is Poincar\'e dual to
$$\sum_{i=1}^{3}r(\gamma_{i})C_{\gamma_{i}}+r(K)C_{K}=-C_{\gamma_{2}}-C_{\gamma_{3}}-\ldots-C_{\gamma_{m+2}}.$$

The cores of the 2-handles form a basis of the 2-chain groups
$C_{2}(W;\mathbb{Z})$; $X_{1}$, $X_{2}$, $\ldots$,  $X_{m+2}$
together  form a basis for the 1-chain groups $C_{1}(W;\mathbb{Z})$,
and the boundary map sends $d_{2}(S_{\gamma_{1}})=X_{2}$,
$d_{2}(S_{\gamma_{2}})=X_{1}-X_{2}$,
$d_{2}(S_{\gamma_{3}})=X_{2}-X_{3}$,
$d_{2}(S_{\gamma_{4}})=X_{3}-X_{4}$, $\ldots$,
$d_{2}(S_{\gamma_{m+2}})=X_{m+1}-X_{m+2}$, $d_{2}(S_{K})=-X_{2}$.
Thus, $h$, the generator of $H_{2}(W;\mathbb{Z})$, can be presented
by $S_{K}+S_{\gamma_{1}}$.

So $\left<c_{1}(J),h\right>=\left<c_{1}(J'),h\right>=0$, and hence $rot(K)=0$.

\begin{remark}
Figure \ref{fig:open_book_twist} and the argument in this section are inspired by \cite{be09}.
\end{remark}

\begin{remark} Suppose $F$ is a page of an open book decomposition supporting $(S^{3}, \xi_{\rm std})$, and $K$ is a homologically nontrivial simple closed curve in $F$, then we can Legendrian realize $K$. To compute the rotation number of $K$, we perform some necessary positive stabilizations to $F$ so that the resulting surface can be obtained from a disk by some positive stabilizations. Note that there exists open book decomposition supporting $(S^3,\xi_{\rm std})$ and cannot be obtained from a disk by positive stabilizations, see \cite{bev10}. On the other hand, the Legendrian knot $K$ does not change in this procedure. So we can factor the monodromy of the open book decomposition or the new one as a product of positive Dehn twists. Therefore we can directly apply the method here to compute the rotation number of $K$.
\end{remark}

\section{\textbf{Proof of Theorem 1.5}}

By the classification of Legendrian representatives of torus knots
in \cite{eh01}, there are $n+2$ Legendrian right  handed trefoil
knots $T(2,3)$ with Thurston-Bennequin invariant $-n$. We order them
by $L_1, L_2,\cdots,L_{n+2}$ such that $rot(L_i)=2i-n-3$.

When we perform Legendrian surgery to $(S^{3},\xi_{\rm std})$ along
$L_{i}$, $i=1,2,\ldots,n+2$, we  obtain a Stein fillable contact
structure $\xi_{i}$ on the three-manifold $S^{3}_{-n-1}(T(2,3))$.
Let $c^{+}(\xi_{i})\in HF^{+}(S^{3}_{n+1}(T(-2,3)))$ and
$\hat{c}(\xi_{i})\in \widehat{HF}(S^{3}_{n+1}(T(-2,3)))$ be the
contact invariants of the contact structure $\xi_{i}$. By
\cite[Theorem 1.5]{os05}, $c^{+}(\xi_{i})$ and
$\widehat{c}(\xi_{i})$ are both nonzero for each $i$.

The Heegaard Floer homology groups of $S^3_n(T(-2,3)$ ($n>6$) are (\cite[Proposition 3.2]{os04})
$$HF^{+}(S^{3}_{n+1}(T(-2,3)), s_{0})=\mathcal{T}^{+}\oplus \mathbb{Z},$$
$$HF^{+}(S^{3}_{n+1}(T(-2,3)), s_{j})=\mathcal{T}^{+}, 1\leq j\leq n, $$ where $s_0,s_1,\cdots,s_n$  denote the $n+1$ Spin$^{c}$ structures of
$S^{3}_{n+1}(T(-2,3))$ and $\mathcal{T}^{+}=\mathbb{Z}[U^{-1}]$, while
$$\widehat{HF}(S^{3}_{n+1}(T(-2,3)), s_{0})=\mathbb{Z}_{(\mathcal{T}^{+})}\oplus\mathbb{Z}\oplus \mathbb{Z},$$
$$\widehat{HF}(S^{3}_{n+1}(T(-2,3)),s_{j})=\mathbb{Z}_{(\mathcal{T}^{+})}, 1\leq j\leq n, $$
where $\mathbb{Z}_{(\mathcal{T}^{+})}\cong\mathbb{Z}$ comes from the kernel of the map $U: \mathcal{T}^{+}\longrightarrow\mathcal{T}^{+}$.

Each contact manifold $(S^{3}_{-n-1}(T(2,3)), \xi_{i})$ bounds a
Stein surface $(X, J_{i})$, and $\langle c_{1}(X, J_{i}), H\rangle =
rot(L_{i})$, where $H$ is the capped Seifert surface of $L_{i}$.
Since $rot(L_{i})$ are pairwise distinct, $(X, J_{i})$
$(i=1,2,\ldots,n+2)$ are pairwise non-isomorphic. By Theorems 2 and
4 in \cite{p04}, $\hat{c}(\xi_{i})$, $i=1,2,\ldots,n+2$, are
distinct primitive elements in $\widehat{HF}(S^{3}_{n+1}(T(-2,3)))$.
Each $\hat{c}(\xi_{i})$ lies in $\widehat{HF}(S^{3}_{n+1}(T(-2,3)),
s_{j})$ for some $j=0,1, \ldots, n$. On the other hand, there are
exactly two distinct primitive elements in
$\mathbb{Z}_{(\mathcal{T}^{+})}$ and there are exactly $n+1$
$\mathbb{Z}_{(\mathcal{T}^{+})}$ summands in
$\widehat{HF}(S^{3}_{n+1}(T(-2,3)))$. Moreover, also by Theorem 4 in
\cite{p04}, for $i_{1}\neq i_{2}$, $\hat{c}(\xi_{i_{1}})$ and
$\hat{c}(\xi_{i_{2}})$ cannot both belong to a
$\mathbb{Z}_{(\mathcal{T}^{+})}$ summand. So at least one of
$\hat{c}(\xi_{i})$, $i=1,2,\ldots,n+2$, does not belong to the $n+1$
$\mathbb{Z}_{(\mathcal{T}^{+})}$ summands in
$\widehat{HF}(S^{3}_{n+1}(T(-2,3)))$. We denote this element as
$\hat{c}(\xi_{i_{0}})$. Therefore, $c^{+}(\xi_{i_{0}})$ does not
belong to the $n+1$ $\mathcal{T}^{+}$-summands in
$HF^{+}(S^{3}_{n+1}(T(-2,3)))$. Hence $c^{+}(\xi_{i_{0}})$ does not vanish in $HF_{\rm red}(S^{3}_{n+1}(T(-2,3)))\cong\mathbb{Z}$.

By Theorem 1.2 in \cite{oss05}, the contact manifold $\xi_{i_{0}}$
does not admit a planar open book decomposition. By Theorem 5.10 in
\cite{o09}, the support genus of $L_{i_{0}}$ is positive.

When $1\leq n\leq 6$, by Theorem 5.9 in \cite{o09}, among the
Legendrian right handed trefoil knots with Thurston Bennequin
invariant $-n$, there must be one whose support genus is positive.

Let $L$ be a Legendrian right handed trefoil with Thurston-Bennequin
invariant $1$, then, by Theorem 1.3, the support genus of
$S_{+}^{n_{1}}S_{-}^{n_{2}}(L)$ is $0$, where $n_{i}\geq 1$
$(i=1,2)$. On the other hand, for $n\geq 2$, the support genus of
$S_{+}^{n}(L)$ and $S_{-}^{n}(L)$ are equal, since they only differ
in the orientation. So the support genus of both $S_{+}^{n}(L)$ and
$S_{-}^{n}(L)$ are $1$.


\begin{thebibliography}{10}

\bibitem[AO01]{AO01}
 Selman~Akbulut and Burak~Ozbagci.
 \newblock Lefschetz fibrations on compact Stein surfaces.
 \newblock {Geometry and Topology}, VOlume 5, 319-334, 2001.

\bibitem[BEV10]{bev10}
 Kenneth~L~Baker, John~B~Etnyre and Jeremy~van~Horn-Morris.
 \newblock Cabling, contact structures and mapping class monoids.
 \newblock arXiv: math.GT/1005.1978, 2010.

\bibitem[BE09]{be09}
 John~A~Baldwin and John~B~Etnyre.
 \newblock A note on the support norm of a contact structure.
 \newblock arXiv: math.GT/0910.5021, 2009.

\bibitem[DG04]{dg04}
 Fan~Ding and Hansjorg~Geiges.
 \newblock A Legendrian surgery presentation of contact 3-manifolds.
 \newblock {\em Mathematical Proceedings of the Cambridge Philosophical Society}, Volume 136, Number 3, 583--598, 2004.

\bibitem[E90]{e90}
 Yakov~M~Eliashberg.
 \newblock Topological characterization of Stein manifolds of dimension $>2$.
 \newblock {\em International Journal of Mathematics}, Volume 1, Number 1, 29--46, 1990.

\bibitem[E05]{E05}
 John~B~Etnyre.
 \newblock Legendrian and transversal knots.
 \newblock {\em Handbook of knot theory}, 105--185, Elsevier B. V., Amsterdam, 2005.

\bibitem[E06]{E06}
 John~B~Etnyre.
 \newblock Lectures on open book decompositions and contact structures.
 \newblock {\em Floer homology, gauge theory, and low-dimensional topology}, 103--141, Clay Math. Proc., 5, Amer. Math. Soc.,
Providence, RI, 2006.

\bibitem[EH01]{eh01}
 John~B~Etnyre and Ko~Honda.
 \newblock Knots and contact geometry. I. Torus knots and the figure eight knot.
 \newblock {\em J. Symplectic Geom}, Volume 1, Number 1, 63--120, 2001.

\bibitem[ENV10]{env10}
 John~B~Etnyre, Lenhard~L~Ng and Vera~Vertesi.
 \newblock Legendrian and transversal twist knots.
 \newblock arXiv: math.GT/1002.2400, 2010.

\bibitem[EO08]{eo08}
 John~B~Etnyre and Burak~Ozbagci.
 \newblock Invariants of contact structures from open books.
 \newblock {\em Transactions of the American Mathematical Society}, Volume 360, Number 6, 3133--3151, 2008.

\bibitem[G02]{g02}
　Emmanuel~Giroux.
 \newblock G\'eom\'etrie de contact: de la dimension trois vers les dimensions sup\'erieures.
 \newblock {\em Proceedings of the
International Congress of Mathematicians}, Vol. II (Beijing, 2002),
405--414, Higher Ed. Press, Beijing, 2002.

\bibitem[LM97]{lm97}
 Paolo~Lisca and Gordana~Matic.
 \newblock Tight contact structures and Seiberg-Witten invariants.
 \newblock {\em Inventiones Mathematicae}, Volume 129, Number 3, 509--525, 1997.

\bibitem[O09]{o09}
 Sinem~Onaran.
 \newblock Invariants of Legendrian knots from open book decompositions.
 \newblock arXiv: math.GT/0905.2238, 2009.

\bibitem[OSz04]{os04}
 Peter~S~Ozsv\'ath and Zolt\'an~Szab\'o.
 \newblock Holomorphic disks and three-manifold invariants: properties and applications.
 \newblock {\em Annals of Mathematics, Second Series}, Volume 159, Number 3, 1159--1245, 2004.

\bibitem[OSz05]{os05}
 Peter~S~Ozsv\'ath and Zolt\'an~Szab\'o.
 \newblock Heegaard Floer homology and contact structures.
 \newblock {\em Duke Mathematical Journal}, Volume 129, Number 1, 39--61, 2005.

\bibitem[OSSz05]{oss05}
 Peter~S~Ozsv\'ath, Andr\'as~Stipsicz and Zolt\'an~Szab\'o.
 \newblock Planar open books and Floer homology.
 \newblock {\em International Mathematics Research Notices}, Volume 2005, Issue 54, 3385--3401, 2005.

\bibitem[P04]{p04}
 Olga~Plamenevskaya.
 \newblock Contact structures with distinct Heegaard Floer invariants.
 \newblock {\em Mathematical Research Letters}, Volume 11, Number 4, 547--561, 2004.

\end{thebibliography}
\end{document}